\documentclass{bioinfo}
\copyrightyear{2015} \pubyear{2015}

\access{Advance Access Publication Date: Day Month Year}
\appnotes{Manuscript Category}

\usepackage{url}

\begin{document}
\firstpage{1}

\newtheorem{theorem}{Theorem}
\newtheorem{corollary}{Corollary}
\newtheorem{lemma}{Lemma}
\newtheorem{proposition}{Proposition}
\newtheorem{definition}{Definition}

\newtheorem{example}{Example}
\newtheorem{remark}{Remark}

\newcommand{\TODO}[1]{\begingroup\color{red}#1\endgroup}

\subtitle{Structural bioinformatics}

\title[Sequence-structure relations]{Sequence-structure relations of biopolymers}
\author[Christopher Barrett \textit{et~al}.]
       {Christopher Barrett\,$^{\text{\sfb 1}}$, Fenix W. Huang\,$^{\text{\sfb 1}}$
         and Christian M. Reidys\,$^{\text{\sfb 1}*}$}
\address{$^{\text{\sf 1}}$Virginia Bioinformatics Institute, 1015 Life Sciences Circle,
Blacksburg, VA, USA}

\corresp{$^\ast$To whom correspondence should be addressed.}

\history{Received on XXXXX; revised on XXXXX; accepted on XXXXX}

\editor{Associate Editor: XXXXXXX}

\abstract{
  {\textbf{Motivation:} DNA data is transcribed into single-stranded RNA, which folds into
    specific molecular structures. In this paper we pose the question to what extent
    sequence- and structure-information correlate. We view this correlation as structural
    semantics of sequence data that allows for a different interpretation than conventional
    sequence alignment. Structural semantics could enable us to identify more general
    embedded ``patterns'' in DNA and RNA sequences.\\
    \textbf{Results:} We compute the partition function of sequences with respect to a fixed
    structure and connect
this computation to the mutual information of a sequence-structure pair for RNA secondary
structures. We present a Boltzmann sampler and obtain the {\it a priori} probability of specific
sequence patterns. We present a detailed analysis for the three
{PDB}-structures, 2JXV (hairpin), 2N3R (3-branch multi-loop) and 1EHZ (tRNA). We localize
specific sequence patterns, contrast the energy spectrum of the Boltzmann sampled
sequences versus those sequences that refold into the same structure and derive a
criterion to identify native structures.
We illustrate that there are multiple sequences in the partition function of a fixed structure,
each having nearly the same mutual information, that are nevertheless poorly aligned. This
indicates the possibility of the existence of relevant patterns embedded in the sequences
that are not discoverable using alignments.}\\
\textbf{Availability:}
The source code is freely available at {http://staff.vbi.vt.edu/fenixh/Sampler.zip}\\
\textbf{Contact:} \href{duckcr@vbi.vt.edu}{duckcr@vbi.vt.edu}\\
\textbf{Supplementary information:}
Supplementary material containing additional data tables are available at \textit{Bioinformatics} online.}

\maketitle

\section{Introduction}

2015 is the 25th year of the human genome project. A recent signature publication \citep{Nature:2015}
is a comprehensive sequence alignment-based analysis of whole genome nucleotide sequence
variation across global human populations. Notwithstanding the importance of this achievement,
there is the possibility of information encoded as patterns in the genome that current methods
cannot discover. 

In this paper we study the information transfer from RNA sequences to RNA structures.
This question is central to the processing of DNA data, specifically the
role of DNA nucleotide sequences being transcribed into RNA, stabilized by
molecular folding. In a plethora of interactions it is this specific configuration and
not the particular sequence of nucleotides (aside from, say small docking areas, where
specific bindings occur) that determines biological functionality. 
We find that here are multiple sequences in the partition function of a fixed structure, each
having nearly the same mutual information with respect to the latter, that are nevertheless
poorly aligned. This indicates the possibility of the existence of relevant patterns embedded
in the sequences that are not discoverable using alignments.

RNA, unlike DNA, is almost always single-stranded\footnote{There are double-stranded RNA
viruses} and all RNA is folded. Here we only consider single-stranded RNA.
An RNA strand has a
backbone made of alternating sugar (ribose) and phosphate groups. Attached to each
sugar is one of four bases--adenine ({\bf A}), uracil ({\bf U}), cytosine ({\bf C}),
or guanine ({\bf G}).
There are various types of RNA: messenger RNA (mRNA), ribosomal RNA
(rRNA), transfer RNA (tRNA) and many others.
Recent transcriptomic and bioinformatic studies suggest the existence of numerous
of so called non-coding RNA, ncRNAs, that is RNA that does not translate into protein
\citep{Eddy:01, Cheng:05}.

\begin{figure}
\begin{center}
\includegraphics[width=1.0\columnwidth]{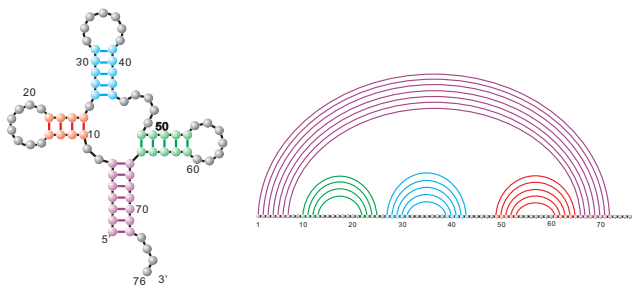}
\end{center}
\caption{\small tRNA: its secondary structure and its diagram presentation. 
}
\label{F:RNAp}
\end{figure}

RNA realizes folded molecular conformations consistent with the Watson-Crick base
as well as the wobble base pairs. In the following we consider RNA secondary
structures, presented as diagrams obtained by drawing the sequence in a
straight line and placing all Watson-Crick and Wobble base pairs as arcs in the upper
half-plane, without any crossing arcs, see Fig.~\ref{F:RNAp}.


DNA information processing refers to replication, transcription and translation.
Additionally, RNA information processing includes replication \citep{Koonin:89}, reverse
transcription (from RNA
to DNA in e.g.~retroviruses) \citep{Temin:70} and a direct translation from DNA to
protein\footnote{in cell-free systems, using extracts from E. coli that contains
ribosomes}\citep{McCarthy:65,Uzawa:02}.

In the following we offer an alternative view of DNA-RNA information processing.
We focus on the information transfer from DNA/RNA sequences to the
folded RNA (after transcription). We speculate that the sequential DNA
information may transcribe into single-stranded RNA in order to allow subsequent
biological processes to interpret DNA data.

DNA data are viewed as sequences of nucleotides. We currently use sequence alignment tools
as a means of arranging the sequences of DNA, RNA, or proteins to identify regions of
similarity that may be a consequence of functional, structural, or evolutionary
relationships between the sequences \citep{Mount:04}. Here we suggest that the transcription
into RNA with the implied self-folding is a way of lifting DNA information to a new and different
level: RNA structures provide sequence semantics.

In order to study this idea we consider the folding of RNA sequences into minimum free
energy (mfe) secondary structures \citep{Waterman:78a}. Pioneered by Waterman more than
three decades ago \citep{Waterman:78aa} and subsequently studied by Schuster {\it et al.}
\citep{Schuster:94} in the context of the RNA toy world
\citep{Schuster:97} there is detailed information about this folding. In particular
we have fairly accurate energy values for computing loop-based mfe \citep{Mathews:99, Mathews:04, Turner:10} that are
employed by the folding algorithms \citep{Zuker:81, Hofacker:94a}. More
work has been done on loop-energy models in \citep{Mathews:04a,Do:06}.
We plan on a more detailed analysis of the framework proposed here in the context of
the MC-model \citep{Parisien:08}.

In \citep{McCaskill:90} McCaskill observed that the dynamic programming routines folding
mfe structures \citep{McCaskill:90} allows one to compute the partition function of all
possible structures for a given sequence. 
The partition function is tantamount to computing the probability space of structures that a
fixed sequence is compatible with. Predictions such as base pairing probabilities are obtained
in \citep{Hofacker:94a, Hofacker:03} and are parallelized in \citep{Fekete:00}. \citep{Tacker:96a, Ding:03}
derives a statistically valid sampling of secondary structures in the Boltzmann ensemble and
calculates the sampling statistics of structural features.  

In view of the above we are led to consider the ``dual'' of McCaskill's partition function, i.e.~the
partition function of all sequences that are compatible with a fixed structure. More generally we
consider the pairing
\begin{equation}
  \varepsilon\colon \mathcal{Q}_4^n \times \mathcal{S}_n \longrightarrow \mathbb{R}^+,
\end{equation}
where $\mathcal{Q}_4^n$ and $\mathcal{S}_n$ denote the space of sequences, $\sigma$, and the
space of secondary structures, $S$, respectively and $\varepsilon(\sigma,S)=e^{-\frac{\eta(
\sigma,S)}{RT}}$ as well as the energy function $\eta(\sigma,S) \in \mathbb{R}$ are discussed in Section~\ref{S:Loops}. 

We show in Section~\ref{S:dis} how $\varepsilon$ allows us to capture the mutual
information between sequences and structures, where the mutual information between 
$x$ and $y$ is given by
$$
I(x,y) = \mathbb{P}(x,y) 
\log \left( \frac{\mathbb{P}(x,y)}{\mathbb{P}(x) \mathbb{P}(y)} \right). 
$$
Here $\mathbb{P}(x,y)$ denotes the joint probability distribution. \TODO{In our case, 
  $\mathbb{P}(\sigma,S) = \epsilon(\sigma,S) / \sum_{\sigma \in \mathcal{Q}_4^n,S\in \mathcal{S}_n}
  \epsilon(\sigma,S)$, $\mathbb{P}(\sigma) = \sum_{S \in \mathcal{S}_n} \mathbb{P}(\sigma,S)$, 
and $\mathbb{P}(S) = \sum_{\sigma\in \mathcal{Q}_4^n} \mathbb{P}(\sigma,S)$.}

In addition, $\varepsilon$ allows us to express folding by considering
$$
\{S \mid \varepsilon(\sigma,S) = \max_{S\in \mathcal{S}_n} \varepsilon(\sigma,S)\},
$$
and inverse folding as to compute
$\{\sigma \mid \varepsilon(\sigma,S)=\max_{S\in \mathcal{S}_n} \varepsilon(\sigma,S)\}$, for
fixed $S$.
Accordingly, the dual to folding is tantamount to computing for fixed $S$
$$
\{\sigma \mid \varepsilon(\sigma,S) = \max_{\sigma\in \mathcal{Q}_4^n} \varepsilon(\sigma,S)\}.
$$
This has direct implications to the
``inverse'' folding of structures. Inverse folding is by construction about the sequence
constraints induced by a fixed structure while avoiding competing configurations. Point
in case: it has been observed in \citep{Busch:06} that starting
with a sequence that is mfe w.r.t.~to a fixed structure, without necessarily folding into it,
constitutes a significantly better initialization than starting with a random sequence.

The paper is organized as follows: we first recall in Section~\ref{S:Loops} 
\TODO{the decomposition of secondary structures
as well as the loop-based thermodynamic model.  }
This in turn facilitates (Sections 2.3 and 2.4) the derivation of the partition function and
Boltzmann sampling. 
In Sections 2.3 and 2.4 we compute $Q(S)$, Boltzmann sampling and the {\it a priori} probability of
sequence patterns.

\begin{methods}
\section{Method}

\subsection{Secondary structures and loop decomposition} \label{S:Loops}

RNA structures can be represented as diagrams where we consider the labels of the sequence
to be placed on the $x$-axis and the Watson-Crick as well as Wobble base pairs drawn as arcs in
the upper half plane see \TODO{Fig.~\ref{F:RNAp}}.
That is, we have a vertex-labeled graph whose vertices are drawn on a
horizontal line labeled by $[n]=\{1, 2, \dots, n\}$, presenting the nucleotides of the
RNA sequence and the linear order of the vertices from left to right indicates the
direction of the backbone from $5'$-end to $3'$-end. Furthermore each vertex can be
paired with at most one other vertex by an arc drawn in the upper half-plane. 
Such an arc, $(i,j)$, represents the base pair between the $i$th and $j$th
nucleotide\footnote{here we assume $j-i>3$ to meet the minimum size requirement of a
hairpin loop.}. 
\TODO{
Two arcs $(i,j)$ and $(r,s)$ are called crossing if and only if  $i<r$ and $i<r<j<s$ holds. 
An RNA structure is called pseudoknot-free, or secondary structure, if it does not contain any 
crossing arcs. Furthermore, the arcs of a secondary structure can be endowed with the partial order: 
$(r,s) \prec (i,j)$ if and only if $i<r<s<j$. 
}

A filtration based on the individual contributions of base pairs of RNA structures was
computed via the Nussinov model \citep{Nussinov:78}. \citep{Waterman:78aa} were the first
bringing energy into the picture, computing the free-energy accurately via loops.
\TODO{A loop in a diagram consists of a sequence of intervals on the backbone $([a_i, b_i])_i$, $1\le i \le k$, 
where $(a_1, b_k)$, $(b_i, a_{i+1})$, for all $1\le i <k$ are base pairs. By construction, each base pair $(i,j)$ is involved 
in exactly two loops: one where $(i,j)$ is maximal respect to $\prec$, and one where $(i,j)$ is not.
Furthermore, there is a distinguished loop, $L_{ex}$, called the exterior loop, where $a_1 = 1$, $b_k = n$ and 
$(a_1, b_k)$ is not a base pair.
}
\TODO{Depending on the number of base pairs, and unpaired bases inside a loop, 
a loop is categorized as hairpin-, containing exactly one base pair and one interval, helix, containing 
 two base pairs and two empty intervals, interior-, containing two non-empty intervals 
and two base pairs, bulge-, containing two base pairs and two intervals, where one of them is empty and 
the other one is not and multi-loops, see Fig.~\ref{F:loops} and Fig.~\ref{F:loop_dia}. 
}

\begin{figure}
\begin{center}
\includegraphics[width=1.0\columnwidth]{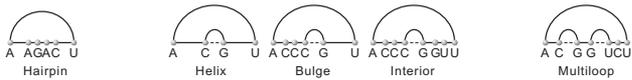}
\end{center}
\caption{\small Hairpin-, helix-, bulge-, interior- and multi-loops in secondary structures.
}
\label{F:loops}
\end{figure}
\TODO{Further 
developments on RNA secondary structure prediction were given by \citep{Zuker:81, Hofacker:94a}.
In particular, accurate thermodynamic energy parameters can be found in
\citep{Mathews:99, Mathews:04, Turner:10, Parisien:08, Deigan:09, Hajdin:13, Lorenz:16}. } 


\begin{figure}
\begin{center}
\includegraphics[width=1.0\columnwidth]{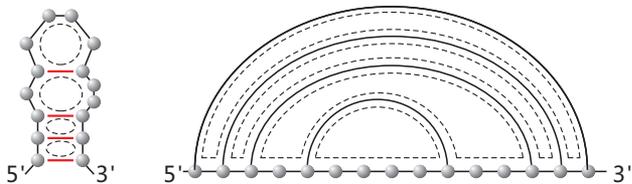}
\end{center}
\caption{\small Loops and their correspondence in a diagram. 
}
\label{F:loop_dia}
\end{figure}

\TODO{In the following, we briefly recall the Turner energy model \citep{Mathews:99, Mathews:04, Turner:10} 
for RNA secondary structures. 
Let $\sigma =(\sigma_1, \sigma_2, \dots, \sigma_n)$ be a sequence,  where 
$\sigma_i \in \{A, U, C, G\}$ for all $1\le i\le n$. }
To an arbitrary loop, $L$, we assign the energy
$\eta(\sigma,L)$, where $\eta(\sigma,L_{ex})=0$ and $\eta(\sigma,L)$ depends
on two factors: its type and the underlying
backbone. Specifically this is the number of bases pairs, the
number of unpaired bases and the particular nucleotides involved. 
The energy of a structure $S$ over an RNA sequence $\sigma$ is then given by
\TODO{the sum of the energies of individual loops i.e., }
\begin{equation}\label{E:xxx}
\eta(\sigma,S) = \sum_{L\in S} \eta(\sigma, L). 
\end{equation}

A hairpin, $L_H$ is a loop having exactly one base pair with a non-empty interval containing $k$
unpaired bases, where $k\ge 3$ due to flexibility constraints imposed by the backbone of the molecule.

In case of $3\le k\le 4$ we call $L$ a tetra-loop, which has a particular energy that depends
on the two nucleotides incident to its unique arc $(\sigma_i, \sigma_{i+k+1})$ as well as the
particular nucleotides corresponding to the unpaired bases of its unique non-empty
interval $(\sigma_{i+1},\dots, \sigma_{i+k})$. 

For any other number of unpaired bases, $k$, the energy calculation depends only on $k$ and not
the particular nucleotide sequence, except of $(\sigma_i,\sigma_{i+k+1})$ and $\sigma_{i+1}$ and 
$\sigma_{i+k}$.  We have 
\begin{equation} \label{E:hp}
\eta(\sigma, L_H) = 
\begin{cases}
\eta_H((\sigma_{i},\sigma_{i+k+1}),\sigma_{i+1},\dots,\sigma_{i+k} ) \quad \text{if $3\le k\le 4$} \\
\eta_{H} ((\sigma_{i},\sigma_{i+k+1}), \sigma_{i+1}, \sigma_{i+k}, k) \quad \text{otherwise}. 
\end{cases}
\end{equation}

An interior, bulge or helix loop, $L_*$, can be represented as two intervals and two base pairs
$L_*= \{[i,r], [s,j], (i,j), (r,s)\}$. 
The energy of $L_*$ is computed as
\begin{equation} \label{E:I}
\eta(\sigma, L_*) = 
\begin{cases}
\eta_{*} ((\sigma_i,\sigma_j),(\sigma_r,\sigma_s)) \quad \text{(helix)}  \\
\eta_{*} ((\sigma_i,\sigma_j),(\sigma_r,\sigma_s),  \sigma_{i+1}, \sigma_{r-1}, \\ \quad
\sigma_{s+1}, \sigma_{j-1} ), k_1) \quad \text {(bulge)} \\
\eta_{*} ((\sigma_i,\sigma_j),(\sigma_r,\sigma_s), \sigma_{i+1}, \sigma_{r-1},  \\ \quad 
\sigma_{s+1}, \sigma_{j-1} ), k_1, k_2) \quad \text {(interior)} 
\end{cases}
\end{equation}
where $k_1= \max \{r-i-1, j-s-1\}$ and $k_2 = \min \{r-i-1, j-s-1\}$.

A multi-loop $L_M$ contains $p$ base pairs and $p$ intervals, some of which being possibly empty, where
$p\ge 3$. $\eta_M(\sigma, L_M)$ is computed by  
\begin{equation} \label{E:M}
\eta_{M} (\sigma, L_M) = \alpha + p \cdot \beta + u \cdot \gamma. 
\end{equation}
Here $\alpha$ is the constant multi-loop penalty, $\beta$ and $\gamma$ are constants
and $u$ is the number of all unpaired bases contained in the respective intervals.

\subsection{The partition function}

\begin{definition}
Let $S$ be a secondary structure over $n$ nucleotides. 
Then the partition function of $S$ is given by
\begin{equation}
Q(S) = \sum_{\sigma\in \mathcal{Q}_4^n} e^{-\frac{\eta(\sigma,S)}{RT}}, 
\end{equation}
where $\eta(\sigma,S)$ is the energy of $S$ on $\sigma$, 
$R$ is the universal gas constant and $T$ is the temperature. 
\end{definition}

\TODO{In analogy to the partition function of a fixed sequence $Q(\sigma)$ \citep{McCaskill:90},
$Q(S)$ can be computed recursively.} Given the structure $S$, we consider an
arbitrary arc $(i, j)$, where $i < j$. Let $S_{i,j}$ denote the substructure of
$S$ over the interval $[i, j]$.
Since $S$ contains no crossing arcs all arcs of $S_{i,j}$ are contained in $[i, j]$,
whence $S_{i,j}$ is well defined. Let
$$
Q (\sigma_i, \sigma_j) = \sum_{\substack{\sigma \in \mathcal{Q}_4^{j-i+1} \\ \sigma|_i= \sigma_i, 
\sigma|_j=\sigma_j}}    e^{-\frac{\eta(\sigma, S_{i,j})}{RT}}. 
$$
Since $S$ has no crossing arcs, \TODO{the interval $[i,j]$ is covered by the arc $(i,j)$, 
i.e.~$(i,j)$ induces a loop $L$ for which $(i,j)$ is maximal. 
Suppose $L$ consist of intervals $[i, p_1], [q_1, p_2], \ldots, [q_k, j]$, where 
$(p_1, q_1) \dots, (p_k, q_k)$ are $L$-arcs different from $(i,j)$.
Removal of $L$ renders substructures covered by $(p_1, q_1) \dots, (p_k, q_k)$.
Considering all combinations of the nucleotides in position $p_i$ and $q_i$, $1\le i \le k$,
we derive the following recursion, see Fig.~\ref{F:rec1}:
} 
\begin{equation}\label{E:xx}
Q (\sigma_i, \sigma_j) = \sum_{\sigma_{p_t}, \sigma_{q_t} \in \mathcal{Q}_4^n}  
e^{-\frac{\eta(\sigma, L)}{RT}} \prod_{t}^k  Q (\sigma_{p_t}, \sigma_{q_t} ). 
\end{equation}

\TODO{The partition function $Q(S)$ is then obtained as the weighted sum of the terms $Q(\sigma_{a_t}, \sigma_{b_t})$, 
where $(a_t,b_t)$, $\forall 1\le t \le k$ are base pairs in the exterior loop $L_{ex}$:
\begin{equation}
Q (S)= \sum_{\sigma_{a_t}, \sigma_{b_t} \in \mathcal{Q}_4^n}  
e^{-\frac{\eta(\sigma, L_{ex})}{RT}} \prod_{t}^k  Q (\sigma_{a_t}, \sigma_{b_t} ). 
\end{equation}
}

{\bf Remark.} The routine of computing $Q(S)$ is similar to the one for finding an optimal
sequence for a given structure in \citep{Busch:06}. Passing to a topological model for RNA
structures \citep{Orland:02,Penner:03,Bon:08,Reidys:11a}, the above 
recursions can be extended to pseudoknotted RNA
structures, i.e.~RNA structures containing crossing arcs. The key here is a general
bijection between maximal arcs and topological boundary components (loops).


\begin{figure}
\begin{center}
\includegraphics[width=1.0\columnwidth]{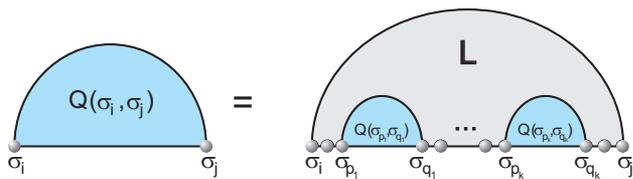}
\end{center}
\caption{\small The recursion for computing the partition function $Q_T(\sigma_i,\sigma_j)$. 
}
\label{F:rec1}
\end{figure}

\subsection{Boltzmann sampling and patterns} 

Having computed the partition function $Q(S)$ as well as the $Q(\sigma_i,\sigma_j)$ terms,
puts us in position to Boltzmann sample sequences for fixed secondary structure
$S$. Here the probability of a sequence $\sigma$ to be sampled is given
by
$$
\mathbb{P}(\sigma|S) = \frac{e^{-\frac{\eta(\sigma,S)}{RT}}}{Q(S)}. 
$$

We build $\sigma$ recursively from top to bottom, starting with the exterior loop,
$L_{ex}$. Suppose $(p_t,q_t)$ are base pairs
contained in $L_{ex}$ and let $u$ denote the number of unpaired bases in $L_{ex}$. 
Since $\eta(\sigma,L_{ex})= 0$, the unpaired nucleotides in $L_{ex}$ are sampled uniformly, 
i.e., with probability $1/4$.  Then 
the probability of the event \TODO{$\sigma_r$ being the nucleotide in position $r\in L_{ex}$}, 
is given by
{\small $$
\mathbb{P}(\sigma_r |S) = \frac{e^{-\frac{\eta(\sigma,L_{ex})}{RT}} \prod_{t=1}^k Q(\sigma_{p_t},\sigma_{q_t}) }
{Q(S)} = \frac{\left(\frac{1}{4}\right) ^u  \prod_{t=1}^k Q(\sigma_{p_t},\sigma_{q_t}) }
{Q(S)}, 
$$
}
where the dependence on $\sigma_r$ of the RHS stems from $\sigma|_r=\sigma_r$ or potentially $p_t=r$
or $q_t=r$.
We continue the process inductively from top to bottom. Suppose we are given a loop $L$ 
with the maximal base pair $(i,j)$. Since any two arcs in $S$ are not crossing, \TODO{any arc
$(i,j)$ is contained in exactly two loops (except for the exterior loop) where $(i,j)$ is the
maximal arc for one and not for the other.} As a result, the nucleotides
$\sigma_i,\sigma_j$ associated with $(i,j)$ are sampled as part of the preceding
loop\footnote{in which $(i,j)$ is not maximal}.
It remains to sample the nucleotides other than $\sigma_i$ and $\sigma_j$ in $L$. 
\TODO{Let $\sigma_r$ be the nucleotides in $L$ and $r\neq i, j$. }
The probability of the event \TODO{$\sigma_r$ being the nucleotide in position $r$}, $r\neq i,j$
is given by
$$
\mathbb{P}(\sigma_r  |S) = \frac{e^{-\frac{\eta(\sigma,L)}{RT}} 
\prod_{t=1}^k Q(\sigma_{p_t}, \sigma_{q_t)}} {Q(\sigma_i, \sigma_j)}.   
$$
Here $(p_t,q_t)$, for $1\le t \le k$, $k\ge 0$ are base pairs contained in $L$, that are different
from $(i,j)$. In particular, $L$ is a hairpin loop in case of $k=0$, an interior-, bulge- or 
a helix-loop in case of $k=1$, and a multi-loop for $k\ge 2$. 

\TODO{By construction, for any arc there is a unique loop
for which the arc is maximal and a unique loop where the arc is not. }
As a result, the probability of a sequence $\sigma$ to be sampled is given by
{\small
\begin{equation*}
\mathbb{P}(\sigma|S)   = \prod_{(i,j)\in S} \frac{e^{-\frac{\eta(\sigma,L(i,j))}{RT}} 
  \prod_{t=1}^k Q(\sigma_{p_t}, \sigma_{q_t})}
       {Q(\sigma_i, \sigma_j)} \cdot \frac{\left(\frac{1}{4}\right) ^u  \prod_{t=1}^k Q(\sigma_{p_t},\sigma_{q_t})}{Q(S)} \\
\end{equation*}
}
In view of eq.~(\ref{E:xxx}) and the fact that the term $Q(\sigma_{p_t,q_t})$ appears exactly once for each
arc $(p_t,q_t)$, we arrive at
{\small
  \begin{equation*}
\prod_{(i,j)\in S}\left(\prod_{t=1}^kQ(\sigma_{p_t},\sigma_{q_t})\right)=\prod_{(i,j)\in S}Q(\sigma_i,\sigma_j).
  \end{equation*}
}
This in turn implies
{\small
  \begin{equation*}    
  \mathbb{P}(\sigma|S)= \frac{\left( \prod_{(i,j)\in S} e^{-\frac{\eta(\sigma,L(i,j))}{RT}}\right) \left(\prod_{(i,j)\in S} Q(\sigma_i,\sigma_j)\right)} 
{Q(S)\prod_{(i,j)\in S} Q(\sigma_i,\sigma_j)} 
=\frac{e^{-\frac{\eta(\sigma,S)}{RT}}}{Q(S)}. 
\end{equation*}
}

The time complexity for computing the partition function of a structure and 
Boltzmann sampling depends solely on the complexity of the energy function,
$\eta(\sigma, L)$. Clearly, there are $O(n)$ loops in the structure and
reviewing eq.~(\ref{E:hp}), eq.~(\ref{E:I}) and eq.~(\ref{E:M}), at most
eight nucleotides are taken into account. From this we can conclude that
the time complexity is $O(n)$, as claimed.

Next, we compute the probability of a given sequence pattern, i.e.~the
subsequence of $\sigma$ over $[i,j]$ being $p_{i,j}$. We shall refer to a
sequence containing $p_{i,j}$ by $\sigma|_{p_{i,j}}$.

The partition function of all sequences $\sigma$ containing $p_{i,j}$ is given by
\begin{equation}\label{E:pattern}
Q(S|p_{i,j}) = \sum_{\sigma |_{p_{i,j}} \in \mathcal{Q}^n_4} 
e^{-\frac{\eta(\sigma, S)}{RT}} 
\end{equation}
and the probability of $p_{i,j}$ is  
$\mathbb{P} (p_{i,j}|S) = \frac{Q(S|p_{i,j})} {Q(S)}$.

We have shown how to compute $Q(S)$ recursively in Section 2.3. It remains to show 
how to compute $Q(S|p_{i,j})$. To do this we use the same routine as for computing 
$Q(S)$, but eliminating any subsequences that are not compatible with $p_{i,j}$.
By construction, for any pattern, this process has the same time complexity as computing 
$Q(S)$. 
\end{methods}

\section{Discussion}\label{S:dis}


Let us begin by discussing the mutual information of sequence-structure pairs. Then we ask
to what extent does a structure determine particular sequence patterns and finally
derive a criterion that differentiates native from random structures. 

The mutual information of a sequence-structure pair can be computed by normalizing $\epsilon$ 
$$
\mathbb{P}(\sigma,S) = \frac{e^{-\frac{\eta(\sigma,S)}{RT}}}{\sum_{\sigma \in  \mathcal{Q}^n_4} \sum_{S\in \mathcal{S}_n}
e^{-\frac{\eta(\sigma,S)}{RT}} },
$$
where $U = \sum_{\sigma \in  \mathcal{Q}^n_4} \sum_{S\in \mathcal{S}_n}
e^{-\frac{\eta(\sigma,S)}{RT}}$ is a constant. 
Then we have 
$$
I(\sigma,S)
=
\left( e^{-\frac{\eta(\sigma,S)}{RT}} \log \frac{(e^{-\frac{\eta(\sigma,S)}{RT}})}{Q(\sigma)Q(S)}\right)/U
  +\left(e^{-\frac{\eta(\sigma,S)}{RT}} {\log U}\right)/U.
$$
Since $U$ is a large constant, we observe that one term of $\mathbb{P}(\sigma,S)$,
namely
$$
e^{-\frac{\eta(\sigma,S)}{RT}} \log \frac{e^{-\frac{\eta(\sigma,S)}{RT}}}
{Q(S)Q(\sigma)}
$$
contributes the most. 
Accordingly, $Q(S)$ and $Q(\sigma)$ allow us to quantify how a probability space of
structures determines a probability space of sequences.

\begin{figure}
\begin{center}
\includegraphics[width=1.0\columnwidth]{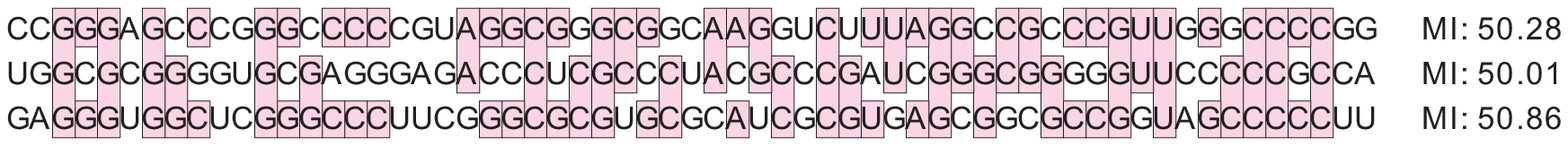}
\end{center}
\caption{\small Three sequences, having nearly the same mutual information with respect to
  the PDB structure 2N3R. The sequences differ pairwise by more than $50$\% of
  their nucleotides which indicates that there is information that cannot be captured by
  conventional sequence alignment. Accordingly BLAST outputs no significant homology
  between the sequences.}
\label{F:MI}
\end{figure}

\begin{figure}[ht]
\begin{center}
\includegraphics[width=1.0\columnwidth]{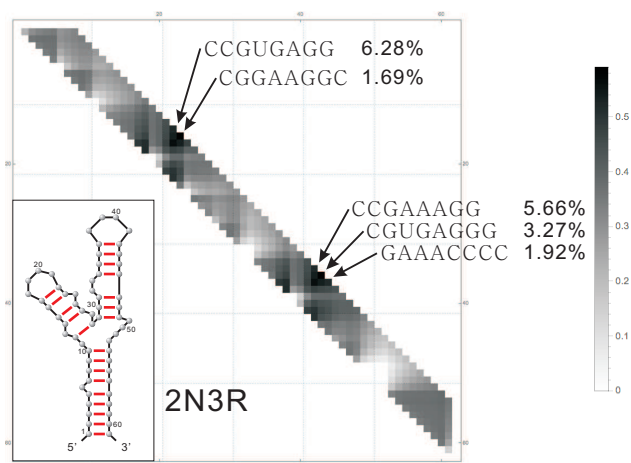}
\end{center}
\caption{\small The secondary structure of 2N3R and its heat-map. We show 
the most frequent patterns for the largest interval having $R_{i,j}>0.52$.  
The sample size is $10^4$. 
}
\label{F:2N3R_h}
\end{figure}

In Fig.~\ref{F:MI} we display three sequences sampled from $Q(S)$ where $S$ is the
PDB-structure 2N3R \citep{PDB:2n3r}, see Fig.~\ref{F:2N3R_h}. All three sequences
have similar mutual information and more than $50\%$ of the nucleotides in the sequences
are pairwise different.


\begin{figure}[ht]
\centerline{\includegraphics[width=1.0\columnwidth]{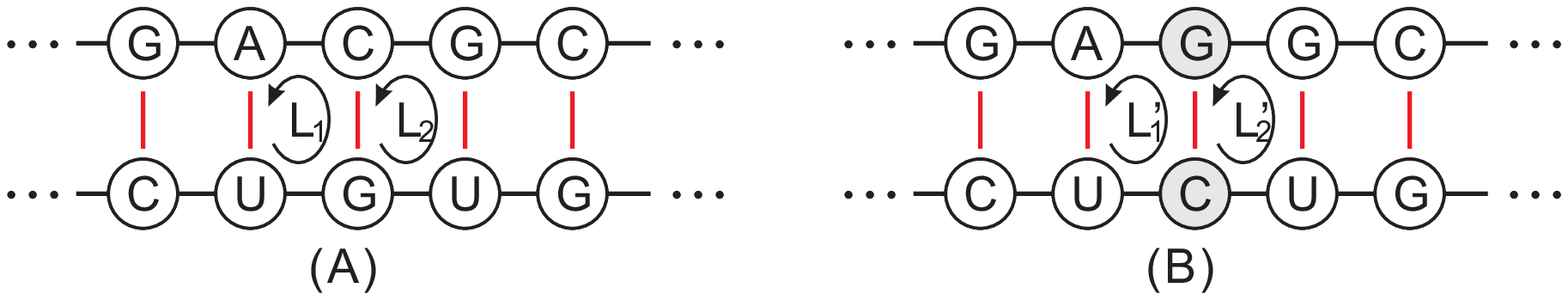}}
\caption{\small
  Isolated replacement of {\bf G}-{\bf C} by {\bf C}-{\bf G}:
  (A) $L_1=({\bf U},{\bf G},{\bf C},{\bf A})$ and
      $L_2=({\bf G},{\bf U},{\bf G},{\bf C})$, 
  (B) replacement induces the new loops: $L_1'=(
           {\bf U},{\bf C},{\bf G},{\bf A})$ and
      $L_2'=({\bf C},{\bf U},{\bf G},{\bf G})$, which changes
      the free energy.  
}
\label{F:GC}
\end{figure}

We point out that replacing a {\bf G}-{\bf C} base pair in a helix by a {\bf C}-{\bf G} base
pair does change the energy, see Fig.~\ref{F:GC}.
This due to the fact that the loops are traversed
in a specific orientation. The isolated replacement of {\bf G}-{\bf C} by {\bf C}-{\bf G} changes
this sequence and hence the energy.

Since the energy model underlying the current analysis does not take non-canonical base pairs
into account, we defer a detailed analysis of the mutual information to a later study where we
use the MC-model \citep{Parisien:08}.

\TODO{Let $p_{i,j}$ be a subsequence on the interval $[i,j]$ with concrete nucleotides}, having probability
$\mathbb{P}(p_{i,j})$. Its Shannon entropy $E_{i,j}$ is given by 
$$
E_{i,j} = - \sum_{\forall p_{i,j}} \mathbb{P}(p_{i,j}) \log_4  \mathbb{P}(p_{i,j}). 
$$
By construction, $0\le E_{i,j} \le (j-i+1)$, where $E_{i,j} = (j-i+1)$ when all $p_{i,j}$
have the same probability, i.e., uniformly distributed, and  $E_{i,j} = 0$ when 
$p_{i,j}$ is completely determined, i.e., $\mathbb{P}(p_{i,j}) = 1$. 
Let $R_{i,j}=1-(E_{i,j}/(j-i+1))$ be the heat of $[i,j]$, i.e.~$R_{i,j}=0$
for random sequences and $R_{i,j}=1$ if there exists only one pattern $p_{i,j}$.
We display the collection of $R_{i,j}$ as a matrix (heat-map), in which we display
$R_{i,j}=0.59$ as black and $R_{i,j}=0$ as white. For a proof of concept,
we restrict ourselves to $R_{i,j}$ for $j-i+1\le 8$. 

\begin{figure}[ht]
\begin{center}
\includegraphics[width=1.0\columnwidth]{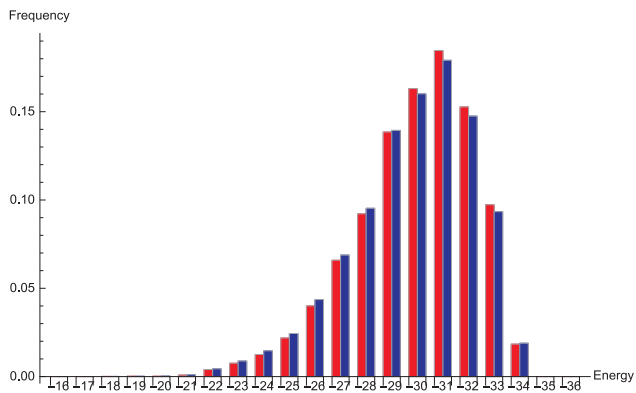}
\end{center}
\caption{\small The energy distribution of the Boltzmann sample for 2JXV.
  We display the frequency of sequences having a particular energy (blue)
  and the frequency of sequences that fold into 2JXV (red).  
}
\label{F:2JXV_e}
\end{figure}

\begin{figure}[ht]
\begin{center}
\includegraphics[width=1.0\columnwidth]{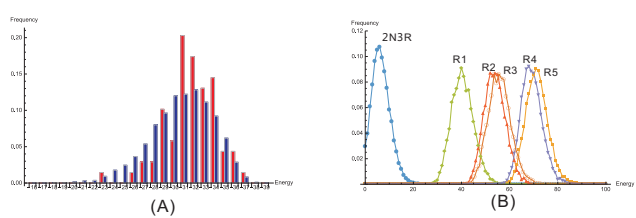}
\end{center}
\caption{\small
(A) The energy distribution of Boltzmann sampled sequences. 
The frequency of sequences having a particular energy level (blue),  
the frequency of sequences folding into 2N3R (red).
(B) $\Delta\eta(\sigma)$-data of 2N3R versus $\Delta\eta(\sigma)$-data of five random
    structures.
}
\label{F:2N3R_e}
\end{figure}

\begin{figure}[ht]
\begin{center}
\includegraphics[width=1.0\columnwidth]{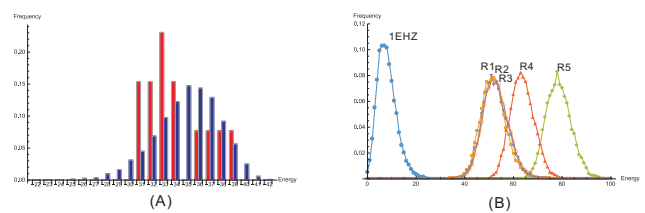}
\end{center}
\caption{\small
(A) The energy distribution of the Boltzmann sampled sequences. 
The frequency of sequences having a particular energy level (blue),  
the frequency of sequences folding into 1EHZ (red).
(B) $\Delta\eta(\sigma)$-data of 1EHZ versus $\Delta\eta(\sigma)$-data of five random
    structures.
}
\label{F:1EHZ_e}
\end{figure}

The heat-maps presented here are obtained by Boltzmann sampling an ensemble of $10^4$
sequences from $Q(S)$. We present the energy distribution of this ensemble in Fig.~\ref{F:2JXV_e},
Fig.~\ref{F:2N3R_e} (A) and Fig.~\ref{F:1EHZ_e} (A)
and in addition the energy spectrum of those sequences that actually fold into $S$ via
the classic folding algorithm using the same energy functions here.  
The Inverse folding rate (IFR),  
$$
\text{IFR} = \frac{\text{\# of sequences folding into $S$}} 
{\text{\# of sampled sequences}}
$$
measures the rate of successful re-folding from that ensemble.

Let $\sigma$ be a sequence from a Boltzmann sample 
w.r.t.~the structure $S$. Let $\bar S$ denote the structure that
$\sigma$ folds to. We consider 
$$
\Delta\eta(\sigma) = |\eta(\sigma,S) -   \eta(\sigma,\bar S)|
$$
and compare the $\Delta G(\sigma,S)$ of several native structures
contained in PDB with those of a several random structures
(obtained by uniformly sampling RNA secondary structures).

\begin{figure}[ht]
\begin{center}
\includegraphics[width=1.0\columnwidth]{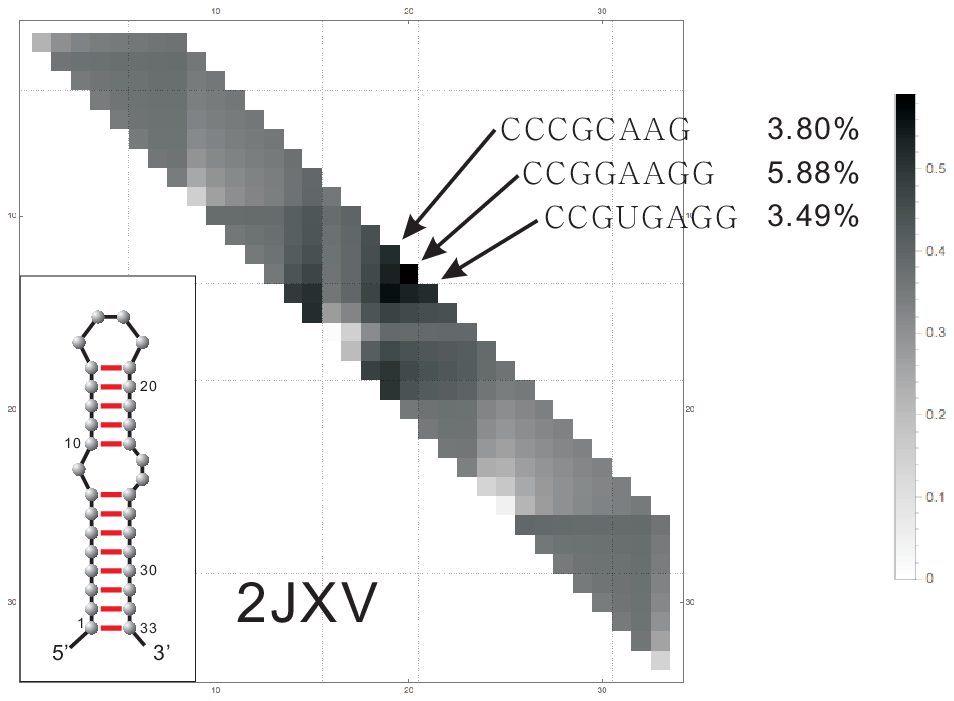}
\end{center}
\caption{\small The secondary structure of 2JXV and its heat-map. We display 
the most frequent sampled patterns for the largest interval having $R_{i,j}>0.52$.  
The sample size is $10^4$. 
}
\label{F:2JXV_h}
\end{figure}

The PDB structure 2JXV \citep{PDB:2jxv} represents a segment of an mRNA, having length $33$. 
The structure exhibits a tetra-loop, an interior loop and 
two stacks of length $8$ and $5$, respectively, see Fig.~\ref{F:2JXV_h}. 
We Boltzmann sample $10^4$ sequences for this structure observing an AU ratio of $18.18\%$, 
while CG ratio is $81.82\%$. The IFR reads $95.16\%$, i.e.~almost all sampled sequences
refold into 2JXV.
The heat-map of 2JXV is given in Fig.~\ref{F:2JXV_h} and we list the most frequent
$10$ patterns of the largest interval having $R_{i,j}>0.52$ in Tab.~ 1 in supplement material 
together with their
{\it a priori} pattern probabilities. We observe that the tetra-loop determines specific patterns.
This finding is not entirely straightforward as the hairpin-loops are the last to be encountered
when Boltzmann sampling. I.e.~they are the most correlated loop-types in the sense that structural
context influences them the most.

The energy distribution of the Boltzmann sample is presented 
in Fig.~\ref{F:2JXV_e} and we observe that the inverse folding solution is not simply the one 
that minimizes the free energy w.r.t.~2JXV with the best energy. 
$\Delta\eta(\sigma)$-data are not displayed here in view of the high IFR.


The PDB structure 2N3R \citep{PDB:2n3r}
consist of $61$ nucleotides and has a $3$-branch multi-loop, two tetra-loops,
interior loops and helixes, see Fig.~\ref{F:2N3R_h}. 
The ratios of AU and CG pairs are $19.67\%$ and $80.33\%$,
respectively, again in a Boltzmann sample of $10^4$ sequences.
The IFR is at $0.69$ quite high, despite the fact that 2N3R is much more complex than 2JXV.
We illustrate the heat-map of 2N3R in Fig.~\ref{F:2N3R_h} and 
list the most frequent $10$
patterns in the largest interval having $R_{i,j}>0.52$ in the Supplemental Materials, Tab. 2
together with their {\it a priori} pattern probabilities computed by eq.~(\ref{E:pattern})

Comparing the sequence segments $[17,24]$ and $[37, 44]$, both of which being tetra-loops with additional
two nucleotides. The $R_{i,j}$ values of these segments are similar, approximately $0.59$, however,
their most frequently sampled patterns appear at different rates. For $[17,24]$ this pattern is
{\tt CGGAAGGC} and it occurs with a Boltzmann sampled frequency of $1.69\%$ and pattern probability
$1.44\%$, while for $[37,44]$ it is {\tt CGUGAGGG} with sampled frequency $3.27\%$ and pattern 
probability $3.24\%$. This makes the point that pattern frequency distributions are strongly
correlated with structural context.

The energy distribution of the Boltzmann sample is given in Fig.~\ref{F:2N3R_e} (A) and we display
the $\Delta\eta(\sigma)$-data in Fig.~\ref{F:2N3R_e} (B) where we contrast the data with
$\Delta\eta(\sigma)$-values obtained from Boltzmann sampling $10^4$ sequences of $5$ random
structures of the same length. We observe that the $\Delta\eta(\sigma)$-values for
2N3R are distinctively lower than those for random structures.
 
\begin{figure}[ht]
\begin{center}
\includegraphics[width=1.0\columnwidth]{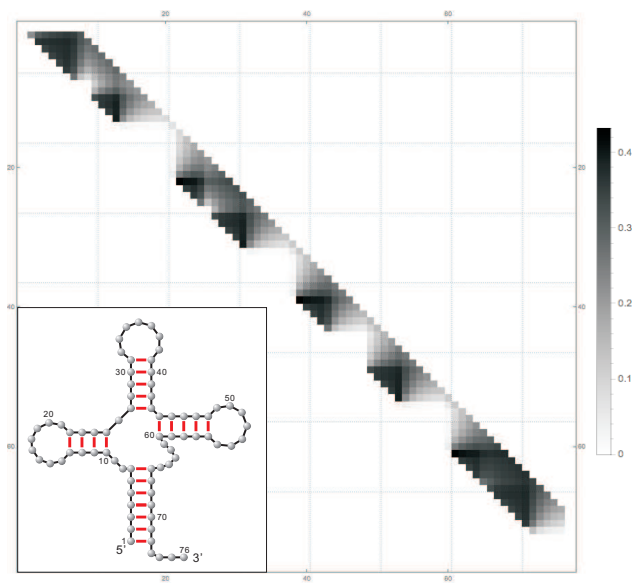}
\end{center}
\caption{\small The secondary structure of 1EHZ and its heat-map. We display
the most frequent sampled pattern for the largest interval having $R_{i,j}>0.52$.  
The sample size is $10^4$.   
}
\label{F:1EHZ_h}
\end{figure}

The PDB structure 1EHZ \citep{PDB:1ehz} is a tRNA over $76$ nucleotides exhibiting a
$4$-branch multi-loop. We display the heat-map of 1EHZ in Fig.~\ref{F:1EHZ_h} 
The IFR is $1.3 \times 10^{-3}$ w.r.t.~our Boltzmann sample of size $10^4$ and we display
the energy distribution of the sampled sequences in Fig.\ref{F:1EHZ_e} (A). Interestingly we
still find many inverse fold solutions by just Boltzmann sampling $Q(S)$ and these
sequences are not concentrated at low free energy values.

In Fig.~\ref{F:1EHZ_e} (B) we display the $\Delta\eta(\sigma)$-data and contrast them with those
obtained by the Boltzmann samples of five random structures. We observe a significant
difference between the $\Delta\eta(\sigma)$-distribution of the 1EHZ sample
and those of the random structures.

The three above examples indicate that sequence-structure correlations can be used to locate
regions where specific embedded patterns arise. Furthermore we observe that studying
$Q(S)$ has direct implications for inverse folding. This is in agreement with
the findings in \citep{Busch:06}, \TODO{but leads to deriving alternative, unbiased
starting sequences for inverse folding.}
Although at present we can only estimate the mutual information, we can conclude
that there are sequences that cannot be aligned but obtain almost identical mutual
information.

\TODO{
We observe that biological relevant sequences exhibit a $\Delta\eta(\sigma)$-signature 
distinctive different from that of random structures. Therefore, the $\Delta\eta(\sigma)$-signature is
capable of distinguishing biological relevant structures from random structures. In
\citep{Miklos:05} the expected free energy and variance of the Boltzmann ensemble of a given
sequence has been employed in order to distinguish biologically functional RNA sequences from
random sequences. 
This result is in terms of the pairing $\varepsilon\colon
\mathcal{Q}_4^n \times \mathcal{S}_n \longrightarrow \mathbb{R}^+$, dual\footnote{the flip side of
  the coin, so to speak} to our approach. Our $\Delta\eta(\sigma)$-signature characterize the naturality
of a fixed structure and \citep{Miklos:05} the naturality of a fixed sequence.
Accordingly, $Q(S)$ augments the analysis of $Q(\sigma)$ in a natural way, capturing the
correlation between RNA sequences and structures.  
}

As a result, sequences carry embedded patterns that cannot be understood by considering
the sequence of nucleotides. At this point we have no concept of what these patterns are and
provided in Section 2.4 a rather conventional notion of ``embedded pattern''. However, even
when considering specific nucleotide patterns in hairpin loops, we observe significant context
dependence on the structure. Other loops affect the energy of the hairpin loop and thus
determine this particular subsequence.
We observe that the embedded patterns can, for certain structures, be quite restricted, possibly
elaborate and are not entirely obvious. In any case, the analysis cannot be reduced to conventional
sequence alignment. The heat-maps introduced here identify the regions
for which only a few select patterns appear and computed the {\it a priori} probabilities of
their occurrence. 

This type of analysis will be carried out for the far more advanced MC-model \citep{Parisien:08},
incorporating non-canonical base pairs, \TODO{SHAPE-directed model for long RNAs \citep{Deigan:09,Hajdin:13}}. 
This will in particular enable us to have a
closer look at the hairpins of the tRNA structure. In addition we believe that this line of
work may enable us to arrive at non-heuristic inverse foldings.

\begin{figure}[ht]
\begin{center}
\includegraphics[width=1.0\columnwidth]{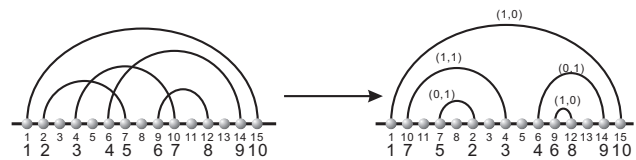}
\end{center}
\caption{\small Disentanglement: by means of permuting the backbone of a pseudoknotted structure one resolves
  all crossings. 
}
\label{F:comm}
\end{figure}

\TODO{
  Folding of RNA secondary structures including pseudoknots is studied in \citep{Rivas:99} by
  extending the dynamic programming paradigm introducing substructures with a gap. The framework generates a
  particular, somewhat subtle class of pseudoknot structures, discussed in detail in \citep{Rivas:00}.
  A specific, multiple context-free grammar (MCFG) for pseudoknotted structures is designed \citep{Rivas:99}, employing
  a vector of nonterminal symbols referencing a substructure with a gap.

  Our results facilitate the Boltzmann sample RNA sequences for pseudoknotted structures. 
  Let $S_{i,j;r,s}$ denote a substructure with a gap where $(i,j)$, $(r,s)$ are base pairs and
  $Q(\sigma_i,\sigma_j;\sigma_r,\sigma_s)$ denote the partition function of $S_{i,j;r,s}$, then one can
  compute $Q(\sigma_i,\sigma_j)$ following the MCFG given by \citep{Rivas:99}. 

  A different approach was presented in \citep{Waterman:93,Penner:03}, where topological RNA structures
  have been introduced \citep{Waterman:93,Penner:03}.
  In difference to \citep{Rivas:99}, which was driven by the dynamic programming paradigm, topological structures stem from
  the intuitive idea to just ``draw'' their arcs on a more complex topological surface in order to
  resolve crossings.
  Random matrix theory \citep{Neumann:47} facilitates the classification and expansion of pseudoknotted structures in terms
  of topological genus \citep{Orland:02,Bon:08} and in \citep{Reidys:11a} a polynomial time, loop-based
  folding algorithm of topological RNA structures was given. The results in this paper are for
  representation purposes formulated in terms of loops. However they were originally developed in the
  topological framework, in which loops become topological boundary components. This means that we can
  extend our framework to pseudoknot structures. The key then is of course to be able to recursively
  compute the novel partition function, i.e.~an unambiguous grammar. Recent results \citep{Huang:16b}
  associate a topological RNA structure with a certain, arc-labeled secondary structure, called
  $\lambda$-structure. 
  The resulting disentanglement gives rise to a context free grammar for RNA pseudoknot structures
  \citep{Huang:16b}\footnote{More precisely, a $\lambda$-structures corresponds one-to-one to a pseudoknotted
  structure together with some additional information, i.e.~a specific permutation of its backbone}. We
  illustrate this correspondence in Fig~\ref{F:comm}. This finding facilitates to extend all our results
  to pseudoknotted structures and offers insight in patterns and inverse folding of more general
  RNA structure classes as well as RNA-RNA
  interaction complexes.
}

As mentioned above, the present analysis is just a first step and discusses embedded patterns
in the sense of subsequent nucleotides. However our framework can deal with any embedded pattern.  
We think a deeper, conceptual analysis has to be undertaken aiming at identifying how a collection
of structures provides sequence semantics. Quite possibly this can be done in the context of
formal languages. We speculate that advancing this may lead to a novel class of embedded pattern
recognition algorithms beyond sequence alignment.

\section{Acknowledgments}
Special thanks to Michael Waterman and Peter Stadler for their input on this manuscript. 
We gratefully acknowledge the help of Kevin Shinpaugh and the computational support team at VBI, 
Rebecca Wattam, Henning Mortveit, Madhav Marathe and Reza Rezazadegan for discussions.
\TODO{Many thanks for the constructive feedback of the anonymous reviewers and their suggestions.}

\bibliographystyle{natbib}

\end{document}